\documentclass[11pt]{article}
\usepackage[margin=1.2in]{geometry}

\usepackage{amsmath}
\usepackage{amsthm}
\usepackage{graphicx}
\usepackage{color}
\usepackage{amsfonts}
\usepackage{amssymb}
\usepackage{amscd}
\usepackage{url}
\usepackage{multirow}%

\newcommand{\re}{\mathbb{R}}

\newcommand{\N}{\mathbb{N}}

\newcommand{\diag}{\mbox{diag}}

\newcommand{\half}{\frac{1}{2}}
\newcommand{\lmd}{\lambda}

\newcommand{\eps}{\epsilon}

\newcommand{\dt}{\delta}
\newcommand{\Dt}{\Delta}
\def\af{\alpha}
\def\bt{\beta}

\def\rank{\mbox{rank}}

\newcommand{\sig}{\sigma}
\newcommand{\Sig}{\Sigma}

\newcommand{\reff}[1]{(\ref{#1})}

\newcommand{\mc}[1]{\mathcal{#1}}

\newcommand{\bdes}{\begin{description}}
\newcommand{\edes}{\end{description}}

\newcommand{\bal}{\begin{align}}
\newcommand{\eal}{\end{align}}

\newcommand{\bnum}{\begin{enumerate}}
\newcommand{\enum}{\end{enumerate}}

\newcommand{\bit}{\begin{itemize}}
\newcommand{\eit}{\end{itemize}}

\newcommand{\bea}{\begin{eqnarray}}
\newcommand{\eea}{\end{eqnarray}}
\newcommand{\be}{\begin{equation}}
\newcommand{\ee}{\end{equation}}

\newcommand{\baray}{\begin{array}}
\newcommand{\earay}{\end{array}}

\newcommand{\bsry}{\begin{subarray}}
\newcommand{\esry}{\end{subarray}}

\newcommand{\bca}{\begin{cases}}
\newcommand{\eca}{\end{cases}}

\newcommand{\bcen}{\begin{center}}
\newcommand{\ecen}{\end{center}}

\newcommand{\bbm}{\begin{bmatrix}}
\newcommand{\ebm}{\end{bmatrix}}

\newcommand{\bmx}{\begin{matrix}}
\newcommand{\emx}{\end{matrix}}

\newcommand{\bpm}{\begin{pmatrix}}
\newcommand{\epm}{\end{pmatrix}}

\newcommand{\btab}{\begin{tabular}}
\newcommand{\etab}{\end{tabular}}

\newtheorem{theorem}{Theorem}[section]

\theoremstyle{definition}

\newtheorem{exm}[theorem]{Example}
\newtheorem{alg}[theorem]{Algorithm}

\newtheorem{remark}[theorem]{Remark}

\begin{document}

\title{Regularization Methods for SDP Relaxations in
Large Scale Polynomial Optimization
\author{Jiawang Nie\footnote{Department of Mathematics,
University of California, 9500 Gilman Drive, La Jolla, CA 92093.
Emails: njw@math.ucsd.edu, liw022@math.ucsd.edu.
The research was partially supported by
NSF grants DMS-0757212 and DMS-0844775.}
\, and \, Li Wang$^*$
}
\date{September 18, 2011}
}

\maketitle

\begin{abstract}
We study how to solve semidefinite programming (SDP) relaxations for
large scale polynomial optimization. When interior-point methods are
used, typically only small or moderately large problems could be solved.
This paper studies regularization methods for
solving polynomial optimization problems. We describe these methods for
semidefinite optimization with block structures, and then apply them
to solve large scale polynomial optimization problems. The
performance is tested on various numerical examples.
By regularization methods, significantly bigger problems could be
solved on a regular computer, which is almost impossible by interior
point methods.
\end{abstract}



\section{Introduction}

Consider the polynomial optimization problem \be \label{pop:S}
\min_{x\in \re^n} \quad f(x) \quad s.t. \quad x \in S \ee where
$f(x)$ is a multivariate polynomial and $S \subseteq \re^n$ is a
semialgebraic set (defined by a boolean combination of polynomial
equalities or inequalities). Recently, there has been much work on
solving \reff{pop:S} by {\it semidefinite programming (SDP)
relaxation} (also called {\it Lasserre's relaxation} in the
literature). The basic idea is approximating nonnegative polynomials
by sum of squares (SOS) type polynomials, which is equivalent to
solving some SDP problems. Thus, the SDP packages (like SDPT3
\cite{sdpt3}, SeDuMi \cite{sedumi}, SDPA\cite{sdpa}) would be
applied to solve polynomial optimization problems.
Typically, SDP relaxation is very successful in solving
\reff{pop:S}, as demonstrated by the pioneer work of Lasserre
\cite{Las01}, Parrilo and Sturmfels \cite{PS01} and many others. However, their
applications are very limited in solving big problems. For instance,
to minimize a general quartic polynomial, it is almost impossible to
solve its SDP relaxation on a regular computer when it has more than $30$ variables.
So far, SDP relaxations for polynomial optimization can only be
solved for small or moderately large problems, which severely
limits their practical applications. Bigger problems would be
solved if sparsity is exploited, like in the work \cite{ND08,WKKM}.
The motivation of this paper is proposing new methods for solving
large scale SDP relaxations arising from general polynomial optimization.

A standard SDP problem is
\be \label{sdp:primal}
\underset{X\in\mc{S}^N}{\min} \quad  C \bullet X
\quad s.t. \quad \mc{A}(X) = b,\, X  \succeq 0.
\ee Here $\mc{S}^N$
denotes the space of $N\times N$ real symmetric matrices, $X\succeq
0$ (resp. $X\succ 0$) means $X$ is positive semidefinite (resp.
definite), and $\bullet$ denotes the standard Frobenius inner
product.
The $C \in \mc{S}^N$ and $b\in \re^m$ are constant, and $\mc{A}: \,
\mc{S}^N \rightarrow \re^m$ is a linear operator. The dual problem
of \reff{sdp:primal} is
\be  \label{sdp:dual}
\max \quad  b^Ty  \,
\quad s.t. \quad \mc{A}^*(y) + Z = C, \, Z \succeq 0.
\ee Here
$\mc{A}^*$ is the adjoint of $\mc{A}$.
An $X$ is optimal for \reff{sdp:primal} and $(y,Z)$ is optimal for
\reff{sdp:dual} if the triple $(X,y,Z)$ satisfies the optimality
condition
\be \label{sdp:optcond} \left. \baray{rl}
\mc{A}(X) & = b \\
\mc{A}^*(y) + Z & = C \\
 X, Z \succeq 0, \,\, XZ & = 0
\earay \right\}.
\ee There is much work on solving SDP by interior
point methods. We refer to \cite{WSV} for theory and algorithms for
SDP. Most of them generate a sequence $\{(X_k,y_k,Z_k)\}$ converging
to an optimal triple.
At each step, a search direction $(\Dt X, \Dt y, \Dt Z)$ needs to be
computed. To compute $\Dt y$, typically an $m\times m$ linear system
needs to be solved. To compute $\Dt X$ and $\Dt Z$, two linear matrix
equations need to be solved. The cost for computing $\Dt y$ is
$\mc{O}(m^3)$. When $m = \mc{O}(N)$, the cost for computing $\Dt y$
is $\mc{O}(N^3)$. In this case, solving SDP is not very expensive if
$N$ is not too big (like less than $1,000$). However, when $m =
\mc{O}(N^2)$, the cost for computing $\Dt y$ would be $\mc{O}(N^6)$,
which is very expensive even for moderately large $N$ (like $500$).
In this case, computing $\Dt y$ is very expensive. It requires
storing a matrix of dimension $m\times m$ in computer and
$\mc{O}(m^3)$ arithmetic operations.

Unfortunately, SDP relaxations arising from polynomial optimization
belong to the bad case that $m=\mc{O}(N^2)$, which is why the SDP
solvers based on interior point methods have difficulty in solving
big polynomial optimization problems (like degree $4$ with $100$ variables).
We explain why this is the case. Let $p(x)$ be a polynomial of
degree $2d$. Then, $p(x)$ is SOS if and only if there exists $X
\succeq 0$ such that $p(x) = [x]_d^T X [x]_d$ (cf. \cite{ParMP}), where
\[
[x]_d^T := [\, 1 \quad  x_1 \quad \cdots \quad x_n \quad x_1^2 \quad
x_1x_2 \quad \cdots \cdots \quad x_1^d \quad x_1^{d-1}x_2 \quad
\cdots \cdots \quad x_n^d \,].
\]
Note the length of $[x]_d$ is
$N=\binom{n+d}{d}$. If we write
\[
 p(x) = \sum_{\af\in\N^n: |\af|
\leq 2d} p_\af x_1^{\af_1}\cdots x_n^{\af_n},
\]
then $p(x)$ being
SOS is equivalent to the existence of a symmetric $N\times N$ matrix
$X$ satisfying
\be  \label{sos:sizes} \baray{rcl}
A_\af \bullet X &=& p_\af \quad \forall \,\af\in\N^n: |\af| \leq 2d, \\
X &\succeq& 0. \earay
\ee Here $A_\af$ are certain constant
symmetric matrices. The number of equalities is $m =
\binom{n+2d}{2d}$. For any fixed $d$,
$m=\mc{O}(n^{2d})=\mc{O}(N^2)$. The size of SDP \reff{sos:sizes} is
huge for moderately large $n$ and $d$. Table~\ref{tab:(N,m)} lists
the size of SDP \reff{sos:sizes} for some typical values of $(n,2d)$.
\begin{table}
\bcen
{\scriptsize
\btab[htb]{|c|c|c|c|c|c|} \hline
n= &  10 & 20 & 30 & 40 & 50  \\
2d = 4&  (66, 1001) & (231, 10626) & (496, 46376) & (861, 135751) & (1326, 316251) \\ \hline
n= &  60 & 70 & 80 & 90 & 100  \\
2d=4 & (1891, 635376) & (2556, 1150626) & (3321, 1929501) & (4186, 3049501) & (5151, 4598126) \\ \hline
n= &  10 & 15 & 20 & 25 & 30  \\
2d = 6& (286, 8008) & (816, 54264) & (1771, 230230) & (3276, 736281) & (5456, 1947792) \\ \hline
n= &  5 & 10 & 15 & 20 & 25  \\
2d = 8& (126, 1287) & (1001, 43758) & (3876, 490314) & (10626,3108105) & (23751,13884156) \\ \hline
n= \, &  5 & 8 & 9 & 10 & 15  \\
2d = 10& (252, 3003) & (1287,43758) & (2002, 92378) & ( 3003, 184756) & (15504,3268760) \\ \hline
\etab}
\caption{A list of sizes of SDP \reff{sos:sizes}.
In each pair (N,m), $N$ is the length of matrix and
$m$ is the number of equality constraints.}
\label{tab:(N,m)}
\ecen
\end{table} 
As we have seen earlier, when interior point methods are
applied to solve \reff{sdp:primal}-\reff{sdp:dual}, at each step we
need to solve a linear system and two matrix equations.
To compute $\Delta y$, we need to store an $m\times m$ matrix and
implement $\mc{O}(n^{6d})$ arithmetic operations. This is very
expensive for even moderately large $n$ and $d$, and hence severely
limits the solvability of SDP relaxations in polynomial
optimization. For instance, on a regular computer, to solve a
general quartic polynomial optimization problem, it is almost impossible to
apply interior point methods when there are more than $30$ variables.

Recently, there has been much work on designing efficient numerical
methods on solving big SDP problems.
{\it Regularization methods} are such a kind of algorithms that
are designed to solve SDP problems whose number of equality constraints $m$
is significantly bigger than the matrix length $N$. We refer to
\cite{MPRW,PRW06,ZST08} for the
work in this area. Their numerical experiments show that these
methods are practical and efficient in solving large scale SDP
problems. In this paper, we study how to apply regularization
methods to solve large scale polynomial optimization problems.

This paper is organized as follows. Section~2 reviews SDP
relaxations in polynomial optimization. Section 3 shows how the
regularization methods work for solving SDP problems whose matrices
have block structures. Section~4 gives numerical experiments in
solving large scale polynomial optimization problems, and Section 5 makes
some discussions about numerical issues.

\bigskip
\noindent {\bf Notations.}  The symbol $\N$ (resp., $\re$) denotes
the set of nonnegative integers (resp., real numbers). For any $t\in
\re$, $\lceil t\rceil$ denotes the smallest integer not smaller than
$t$. For $x \in \re^n$, $x_i$ denotes the $i$-th component of $x$,
that is, $x=(x_1,\ldots,x_n)$. The $\mathbb{S}^{n-1}$ denotes the
$n-1$ dimensional unit sphere $\{x\in\re^n: x_1^2+\cdots+x_n^2=1\}$.
For $\af \in \N^n$, denote $|\af| = \af_1 + \cdots + \af_n$. The
symbol $\N_{\leq k}$ denotes the set $\{\af\in\N^n: |\af| \leq k
\}$, and $\N_{k}$ denotes $\{\af\in\N^n: |\af| = k \}$. For each
$i$, $e_i$ denotes the $i$-th standard unit vector. The $\mathbf{1}$
denotes a vector of all ones.
For $x \in \re^n$ and $\af \in \N^n$, $x^\af$ denotes
$x_1^{\af_1}\cdots x_n^{\af_n}$.
For a finite set $T$, $|T|$ denotes its cardinality. For a matrix
$A$, $A^T$ denotes its transpose. The $I_N$ denotes the $N\times N$
identity matrix, and $\mathcal{S}_+^N$ denotes the cone of symmetric
positive semidefinite $N\times N$ matrices.
For any vector $u\in \re^N$, $\| u \|_2 = \sqrt{u^Tu}$ denotes the
standard Euclidean norm.

\section{SDP relaxations for polynomial optimization}

This section reviews constructions of SDP relaxations for
polynomial optimization problems of three different types:
unconstrained polynomial
optimization, homogeneous polynomial optimization, and constrained
polynomial optimization.

\subsection{Unconstrained polynomial optimization}

\setcounter{equation}{0} \label{sec:unopt}

Consider the optimization problem
\be \label{uc:minf(x)}
f_{min}^{uc}:= \min_{x\in \re^n} \quad f(x)
\ee where $f(x)$ is a
polynomial of degree $2d$, and $f_{min}^{uc}$ denotes the global
minimum of $f(x)$ over $\re^n$.
A standard SOS relaxation for \reff{uc:minf(x)} (cf. \cite{PS01,ParMP}) is
\be \label{opt:stdsos}
f_{sos}^{uc}:=
\max \quad \gamma  \quad
s.t. \quad f(x) - \gamma \mbox{ is SOS}.
\ee Obviously the above
optimal value $f_{sos}^{uc}$ satisfies the relation $f_{sos}^{uc}
\leq f_{min}^{uc}$. Though it is possible that $f_{sos}^{uc} <
f_{min}^{uc}$, it was observed in \cite{PS01,ParMP} that
\reff{opt:stdsos}  works very well in practice.
In the following, we show how to transform \reff{opt:stdsos} into a
standard SDP.

Denote $\mathbb{U}_{2d}^n=\{\af\in \N^n: 0 < |\af| \leq 2d\}$ and write
\[
f(x) = f_0+\sum_{\af\in \mathbb{U}_{2d}^n } f_\af x^\af,
\]
then $f(x)-\gamma$ is SOS if and only if there exists
$X\in\mc{S}^{\binom{n+d}{d}}$ satisfying
\be \label{sos:main}
f(x)-\gamma = [x]_d^TX[x]_d = X \bullet ([x]_d[x]_d^T), \quad X
\succeq 0.
\ee Note that $\binom{n+d}{d}$ is the length of $[x]_d$. Let
$b=(f_\af)_{\af\in\mathbb{U}_{2d}^n}$,  whose dimension is $m
=\binom{n+2d}{2d}-1$. Define 0/1 constant symmetric matrices $C$ and
$A_\af$ such that
\be \label{def:ucpop-A} [x]_d[x]_d^T = C +
\sum_{\af\in\N^n: 0< |\af| \leq 2d} A_\af x^\af.
\ee
Then, (\ref{sos:main}) can be expressed as follows:
\be \label{sos:linear}
f(x)-\gamma = C \bullet X +\sum_{\af\in\N^n: 0< |\af| \leq 2d}
(A_\af \bullet X )x^\af, \quad X \succeq 0.
\ee So, $\gamma$ is
feasible for \reff{opt:stdsos} if and only if
there is a symmetric matrix $X$ satisfying
\[
\baray{rcl}
C \bullet X + \gamma  & = & f_0, \\
A_\af \bullet X & = & f_\af \quad \forall \, \af\in\mathbb{U}_{2d}^n, \\
X & \succeq & 0. \earay
\]
Define a linear operator $ \mc{A}(X) = (A_\af \bullet
X)_{\af\in\mathbb{U}_{2d}^n}. $
Then, up to a constant, SOS relaxation
\reff{opt:stdsos} is equivalent to the SDP problem
\be \label{sdp:stdsos}
f_{sdp}^{uc} :=
\min \quad C \bullet X \quad
s.t. \quad \mc{A}(X) = b, \, X  \succeq 0.
\ee
The dual optimization of the above is
\be  \label{sdp:momsos}
\max \quad b^Ty  \quad
s.t. \quad \mc{A}^*(y) + Z =C, \, Z \succeq 0.
\ee Here
$\mc{A}^*(y) = \sum_{\af \in \mathbb{U}_{2d}^n } y_\af A_\af$.
Clearly, it holds that $f_{sos}^{uc} = -f_{sdp}^{uc}+f_0$.

Suppose $(X^*,y^*,Z^*)$ is an optimal triple for
\reff{sdp:stdsos}-\reff{sdp:momsos}. Then $-f_{sdp}^{uc}+f_0$ is a
lower bound of the minimum $f_{min}^{uc}$. As is well known, if $Z^*$ has
rank one, then $f_{min}^{uc}=f_{sos}^{uc}$ and a global minimizer
for \reff{uc:minf(x)} can be obtained easily. This can be
illustrated as follows. When $\rank(Z^*)=1$, the constraint in
\reff{sdp:momsos} implies $Z^*=[x^*]_d[x^*]_d^T$ for some $x^*\in
\re^n$, and hence $y^*= -[x^*]_{2d}$. Then, for any $x\in \re^n$,
\[
-f(x^*) = -f_0 + b^T y^* \geq -f_0 + \sum_{\af \in
\mathbb{U}_{2d}^n} -b_\af x^{\af} = -f(x).
\]
In the above, we have used the optimality of $y^*$ and the fact that
$Z=[x]_d[x]_d^T$ is always feasible for \reff{sdp:momsos}. So $x^*$
is a global minimizer.

When $\rank(Z^*)>1$, several global minimizers for \reff{uc:minf(x)}
could be obtained if the {\it flat extension condition (FEC)} holds.
We refer to Curto and Fialkow \cite{CF05} for FEC, and Henrion and
Lasserre \cite{HL05} for a numerical method on how to get global
minimizers when FEC holds.
Typically, FEC fails if the SDP relaxation is not exact.

\subsection{Homogeneous polynomial optimization}
  \label{sec:hmgopt}
Consider the homogeneous polynomial optimization problem
\be
\label{hmgpop:sphere}
 f_{min}^{hmg}:=
\underset{x\in\re^n}{\min} \quad f(x) \quad
s.t. \quad \|x\|_2 = 1,
\ee
where $f(x)$ is a form (homogeneous
polynomial). Assume its degree
$\deg(f)=2d$ is even. An interesting application of \reff{hmgpop:sphere} is computing
stability number of graphs \cite{dKP02}. This will also be shown in
Section 4.2.

A standard SOS relaxation for \reff{hmgpop:sphere} is
\be \label{hmgpop:sos}
 f_{sos}^{hmg}:=
\max \quad \gamma \quad
s.t. \quad f(x) - \gamma \cdot (x^Tx)^d \quad \mbox{ is SOS}.
\ee Let $[x^d]$ be the vector of monomials of degree $d$ ordered
lexicographically, i.e.,
\[
[x^d]^T = \bbm x_1^d & x_1^{d-1}x_2 & x_1^{d-1}x_3 & \cdots &
x_{n-1}x_n^{d-1} & x_n^d \ebm.
\]
Denote $\N_{d} = \{\af\in \N^n: |\af| = d\}$. For each $\af \in
\N_d$, define $D_\af = \frac{(|\af|)!}{\af_1! \cdots \af_n!}$. Let
$D=\diag(D_\af)$ be a diagonal matrix. Then, it holds the relation
\be \label{def:Daf} (x^Tx)^d = \sum_{\af \in \N_{d} } D_\af
x_1^{2\af_1}\cdots x_n^{2\af_n} = [x^d]^TD[x^d] =
([x^d][x^d]^T)\bullet D. \ee Define $0/1$ matrices $H_\af$ such that
\be  \label{def:H_af} [x^d][x^d]^T = \sum_{\af \in \N_{2d} }  H_\af
x_1^{\af_1}\cdots x_n^{\af_n}. \ee Then, $f(x)-\gamma \cdot
(x^Tx)^d$ being SOS is equivalent to the existence of $X$ satisfying
\[
\baray{rl}
f_\af - \gamma D_{\af/2} &= H_\af \bullet X \quad \forall\, \af \in 2\N_{d},  \\
f_\af  & = H_\af \bullet X \quad \forall \, \af \in \N_{2d}\backslash 2\N_{d}, \\
X& \succeq 0, \earay
\]
where $f_\af$ is the coefficient of $x^\af$ in $f(x)$.
Letting $\af = 2de_1$ in the above, we get $ \gamma = f_{2de_1} -
H_{2de_1} \bullet X. $ Denote $\mathbb{H}_{2d}^n = \N_{2d}\backslash
\{2de_1\}$, and set $r=(r_\af)_{\af\in \mathbb{H}_{2d}^n}$ as \be
\label{hmg-def:r} r_\af = \bca
D_{\af/2}  & \text{ if } \af \in 2\N_{d}\backslash\{2de_1\}, \\
0 & \text{ if } \af \in \N_{2d}\backslash\{2\N_{d}\}. \eca \ee
Define matrices $C, A_\af$ and scalars $b_\af$ as
\be
\label{matdef:Homg} C= H_{2de_1}, \quad A_\af = H_\af - r_\af
H_{2de_1}, \quad b_\af = f_\af - r_\af f_{2de_1},   \quad
\af\in \mathbb{H}_{2d}^n.
\ee Let $b=(b_\af)_{\af\in \mathbb{H}_{2d}^n}$.
Define linear operators $\mc{H}, \mc{A}:
\mc{S}^{\N_{d}} \rightarrow \re^{\mathbb{H}_{2d}^n}$ as
\[
\mc{H}(X) = (H_\af \bullet X)_{\af\in \mathbb{H}_{2d}^n}, \quad
\mc{A}(X) = (A_\af \bullet X)_{\af\in \mathbb{H}_{2d}^n}.
\]
Then, SOS relaxation \reff{hmgpop:sos} is equivalent to the SDP problem
\be \label{hmgsdp:prime}
f_{sdp}^{hmg}:=
\min \quad C \bullet X \quad
s.t. \quad \mc{A}(X) = b, \, X \succeq 0.
\ee The dual problem of \reff{hmgsdp:prime} is
\be \label{hmgsdp:dual}
\max \quad b^Ty \quad
s.t. \quad  \mc{A}^*(y) +Z = C, \, Z \succeq 0.
\ee In the above,
$\mc{A}^*(y) = \sum_{\af\in \mathbb{H}_{2d}^n} y_\af A_\af$.
Clearly, $f_{sos}^{hmg} = -f_{sdp}^{hmg}+f_{2de_1}$.

Let $(X^*,y^*,Z^*)$ be an optimal triple for
\reff{hmgsdp:prime}-\reff{hmgsdp:dual}. Then
$-f_{sdp}^{hmg}+f_{2de_1}$ is a lower bound of the minimum $f_{min}^{hmg}$.
%
We could also get global
minimizers from $Z^*$ when FEC holds. Note that \reff{def:Daf} and
\reff{matdef:Homg} imply \be  \label{A.D=Z.D=0} A_\af \bullet D = 0
\quad \forall \af \in \mathbb{H}_{2d}^n  \quad \text{ and } \quad
Z^* \bullet D = H_{2de_1} \bullet D = 1. \ee So $Z^*(de_i,de_i)$
($Z$ is indexed by integer vectors in $\N^n$) can not vanish for
every $i$, because otherwise we would get $Z^*=0$ contradicting
\reff{A.D=Z.D=0}. Up to a permutation of $x$, we can assume
$Z^*(de_1,de_1)\ne 0$, and normalize $Z^*$ as
$\widehat{Z^*}=Z^*/Z^*(de_1,de_1)$. Then $\widehat{Z^*}$ would be
thought of as a moment matrix of order $d$ in $n-1$ variables (see
\cite{CF05}). If $\widehat{Z^*}$ satisfies FEC, using the method in
\cite{HL05}, we can get $v^{(1)},\ldots,v^{(r)}\in\re^{n-1}$ such
that
\[
\widehat{Z^*} = \lmd_1 [v^{(1)}]_d[v^{(1)}]_d^T+\cdots
+\lmd_r[v^{(r)}]_d[v^{(r)}]_d^T
\]
for some scalars $\lmd_i > 0$. Setting $x^{(i)} =
\big(1+\|v^{(i)}\|_2^2\big)^{-1/2} \bbm 1 \\ v^{(i)} \ebm \in
\mathbb{S}^{n-1}$, we get
\[
Z^*  = \nu_1 [(x^{(1)})^d] [(x^{(1)})^d]^T + \cdots + \nu_r
[(x^{(r)})^d] [(x^{(r)})^d]^T
\]
for some scalars $\nu_i > 0$. Then, the relations \reff{def:Daf} and
\reff{A.D=Z.D=0} imply $\nu_1+\cdots +\nu_r=1$. Since every
$Z^{(i)}=[(x^{(i)})^d] [(x^{(i)})^d]^T$ is feasible for
\reff{hmgsdp:dual}, the optimality of $Z^*$ implies every $x^{(i)}$
is a global minimizer of \reff{hmgpop:sphere}.

%
%

\subsection{Constrained polynomial optimization}
\label{sec:conopt}

Consider the general polynomial optimization problem \be
\label{conpop:gi>=0} \baray{rl} f_{min}^{con}:=
\underset{x\in \re^n}{\min} & f(x) \\
s.t. & g_1(x)\geq 0, \ldots, g_\ell(x) \geq 0, \earay \ee where
$f(x)$ and $g_1(x),\ldots,g_\ell(x)$ are all polynomials in $x$ and
their degrees are no greater than $2d$.
Problem \reff{conpop:gi>=0} is NP-hard, even when $f(x)$ is
quadratic and the feasible set is a simplex. Lasserre's relaxation
\cite{Las01} is a typical approach for solving \reff{conpop:gi>=0}.
The $d$-th Lasserre's relaxation ($d$ is also called the relaxation
order) for \reff{conpop:gi>=0} is
\be  \label{consos:deg2d}
\baray{rl}
f_{sos}^{con}:=
\max& \gamma \\
s.t. & f(x) - \gamma = \sig_0(x) + g_1(x)\sig_1(x) + \cdots + g_\ell(x)\sig_\ell(x), \\
& \sig_0(x),\sig_1(x),\ldots,  \sig_\ell(x) \mbox{ are SOS}, \\
& \deg(\sig_0),\deg(\sig_1g_1),\ldots,  \deg(\sig_\ell g_\ell) \leq 2d.
\earay
\ee
Let $N(k) = \binom{n+k}{k}$,
$d_i=\lceil\deg(g_i)/2\rceil$ and $g_0(x)=1$. Then, $\gamma$ is
feasible for \reff{consos:deg2d} if and only if there exists
$X^{(i)} \in\mc{S}^{N(d-d_i)}\,(i=0,1,\ldots,\ell)$ such that
\[
\baray{c} f(x)-\gamma = \overset{\ell}{\underset{i=0}{\sum}} g_i(x)
[x]_{d-d_i}^TX^{(i)}[x]_{d-d_i}
= \overset{\ell}{\underset{i=0}{\sum}} X^{(i)} \bullet (g_i(x)[x]_{d-d_i}[x]_{d-d_i}^T), \\
\, X^{(0)} \succeq 0, X^{(1)} \succeq 0, \ldots, X^{(\ell)} \succeq
0. \earay
\]
Define constant symmetric matrices $A_\af^{(0)}, A_\af^{(1)},
\ldots, A_\af^{(\ell)}$ such that \be \label{consos:diagA}
g_i(x)[x]_{d-d_i}[x]_{d-d_i}^T = \sum_{\af\in\N^n: |\af| \leq 2d}
A_\af^{(i)} x^\af, \quad i=0,1,\ldots,\ell. \ee Denote $A_\af =
(A_\af^{(0)}, A_\af^{(1)}, \ldots, A_\af^{(\ell)})$, $X = (X^{(0)},
X^{(1)}, \ldots, X^{(\ell)})$, and define a cone of products
\[
\mc{K} := \mc{S}_+^{N(d-d_0)} \times \mc{S}_+^{N(d-d_1)} \times
\cdots \times \mc{S}_+^{N(d-d_\ell)}.
\]
Recall that $\mathbb{U}_{2d}^n=\{\af\in \N^n: 0 < |\af| \leq 2d\}$.
If $ f(x) =f_0+ \sum_{\af\in\mathbb{U}_{2d}^n} f_\af x^\af$, then
$\gamma$ is feasible for \reff{consos:deg2d} if and only if there
exists $X$ satisfying
\[
\baray{rcl}
A_0 \bullet X + \gamma  & = & f_0, \\
A_\af \bullet X & = & f_\af \quad \forall \, \af\in\mathbb{U}_{2d}^n, \\
X & \in & \mc{K}. \earay
\]
Now define $\mc{A},b,C$ as
\be  \label{consos:data-def} \mc{A}(X) =
(A_\af \bullet X)_{\af\in\mathbb{U}_{2d}^n}, \quad
b=(f_\af)_{\af\in\mathbb{U}_{2d}^n}, \quad C=A_0.
\ee The vector $b$
has dimension $m =N(2d)-1$. Then, up to a constant,
\reff{consos:deg2d} is equivalent
to the SDP problem
\be  \label{prmsdp:consos}
f_{sdp}^{con}:=
\min \quad C \bullet X \quad
s.t. \quad \mc{A}(X) = b, \, X \in \mc{K}.
\ee
Its dual optimization is
\be  \label{dualsdp:consos}
\max \quad  b^Ty  \quad
s.t. \quad \mc{A}^*(y) + Z = C,\, Z \in \mc{K}.
\ee
The adjoint $\mc{A}^*(y)$ is defined as
\[
\mc{A}^*(y) = \sum_{\af\in\mathbb{U}_{2d}^n} y_\af \, \diag \left(
A_\af^{(0)},  A_\af^{(1)}, \ldots,   A_\af^{(\ell)} \right).
\]
Note the relation $f_{sos}^{con} = -f_{sdp}^{con}+f_0 \leq
f_{min}^{con}$.

Suppose $(X^*,y^*,Z^*)$ is an optimal triple for
\reff{prmsdp:consos}-\reff{dualsdp:consos}. Then $-f_{sdp}^{con} +
f_0$ is a lower bound of the minimum $f_{min}^{con}$. The
information for minimizers could be obtained from $Z^*$. Note
$Z^*=(Z_0^*,Z_1^*,\ldots,Z_\ell^*)$. Since $Z^*\in \mc{K}$, every
$Z_i^* \succeq 0$. If $Z_0^*$ satisfies FEC, one or several global
minimizers can be obtained (cf. \cite{HL05}).

\section{Regularization methods }
\setcounter{equation}{0} \label{sec:regalg}

This section describes regularization methods for solving
semidefinite optimization problems having block diagonal structures.
They are natural generalizations of regularization methods
introduced in \cite{MPRW,PRW06,ZST08} for solving standard SDP
problems of a single block structure.


Let $\mc{K}$ be a cross product of several semidefinite cones
\[
\mc{K} = \mc{S}^{N_1}_+ \times \cdots \times \mc{S}^{N_\ell}_+.
\]
It belongs to the space $ \mc{M} = \mc{S}^{N_1} \times \cdots \times
\mc{S}^{N_\ell}. $ Each $X \in \mc{M}$ is a tuple
$X=(X_1,\ldots,X_\ell)$ with every $X_i\in \mc{S}^{N_i}$. So, $X$
could also be thought of as a symmetric block diagonal matrix, and
$X\in\mc{K}$ if and only if its every block $X_i \succeq 0$. The
notation $X \succeq_{\mc{K}} 0$ (resp. $X \succ_{\mc{K}} 0$) means
every block of $X$ is positive semidefinite (resp. definite). For
$X=(X_1,\ldots,X_\ell) \in \mc{M}$ and $Y=(Y_1,\ldots,Y_\ell) \in
\mc{M}$, define their inner product as $ X \bullet Y =  X_1 \bullet
Y_1 + \cdots + X_\ell \bullet Y_\ell. $ Denote by $\| \cdot \|$ the
norm in $\mc{M}$ induced by this inner product.
Note $\mc{K}$ is a self-dual cone, that is,
\[
\mc{K}^*:= \{ Y\in \mc{M}: Y \bullet X \geq 0 \quad \forall X \in
\mc{K} \} = \mc{K}.
\]

For a symmetric matrix $W$, denote by $(W)_+$ (resp. $(W)_{-}$) the
projection of $W$ into the positive (resp. negative) semidefinite
cone, that is, if $W$ has spectral decomposition
\[
W = \sum_{\lmd_i>0} \lmd_i u_iu_i^T + \sum_{\lmd_i< 0} \lmd_i
u_iu_i^T,
\]
then $(W)_+$ and $(W)_{-}$ are defined as
\[
(W)_+ = \sum_{\lmd_i>0} \lmd_i u_iu_i^T, \quad (W)_{-} =
\sum_{\lmd_i< 0} \lmd_i u_iu_i^T.
\]
For $X=(X_1,\ldots,X_\ell) \in \mc{M}$, its projections into
$\mc{K}$ and $-\mc{K}$ are given by
\[
(X)_{\mc{K}} =((X_1)_+,\ldots,(X_\ell)_+), \quad (X)_{-\mc{K}}
=((X_1)_{-},\ldots,(X_\ell)_{-}).
\]

A general conic semidefinite optimization problem is
\be \label{Ksdp:prm}
\min \quad  C \bullet X \quad
s.t. \quad \mc{A}(X) = b, \,  X  \in \mc{K}.
\ee Here $C \in
\mc{M}$, $b\in\re^m$, and $\mc{A}: \mc{M} \to \re^m$ is a linear
operator. The dual of \reff{Ksdp:prm} is
\be \label{Ksdp:dual}
\max \quad  b^Ty \quad
s.t. \quad \mc{A}^*(y) + Z =C, \, Z \in \mc{K}.
\ee SDP
relaxations for constrained polynomial optimization often have block
diagonal structure, e.g., \reff{prmsdp:consos}.

There are two typical regularizations for standard SDP problems: {\it
Moreau-Yosida regularization} for the primal \reff{sdp:primal} and
{\it Augmented Lagrangian regularization} for the dual
\reff{sdp:dual}. They would be naturally generalized to conic
semidefinite optimization problem \reff{Ksdp:prm} and its dual
\reff{Ksdp:dual}. The Moreau-Yosida regularization for
\reff{Ksdp:prm} is \be \label{sdp:MYreg} \baray{rl}
\underset{X,Y \in \mc{M} }{\rm min} &  C \bullet X + \frac{1}{2\sigma} \|X-Y\|^2 \\
s.t. & \mc{A}(X) = b,\, X  \in \mc{K}. \earay \ee Obviously
\reff{sdp:MYreg} is equivalent to \reff{Ksdp:prm}, because for each
fixed feasible $X \in \mc{M}$ the optimal $Y\in\mc{M}$ in
\reff{sdp:MYreg} is equal to $X$. The Augmented Lagrangian
regularization for \reff{Ksdp:dual} is \be \label{sdp:augLag:1st}
\baray{rl} \underset{y\in\re^m, Z\in \mc{M}}{\rm max} &
b^Ty - (Z+\mc{A}^*(y)-C)\bullet Y - \frac{\sigma}{2} \|Z+\mc{A}^*(y) -C\|^2 \\
s.t. &  Z \in \mc{K}. \earay \ee When $\mc{K}=\mc{S}_+^N$ is a
single product, it can be shown (cf. \cite[Section~2]{MPRW}) that for
every fixed $Y$, \reff{sdp:augLag:1st} is the dual optimization
problem of
\[
\min_{X\in\mc{M} } \quad  C \bullet X + \frac{1}{2\sigma} \|X-Y\|^2
-y^T( \mc{A}(X)-b) - Z \bullet X.
\]
By fixing $y\in\re^m$ and optimizing over $Z\succeq 0$, Malick,
Povh, Rendl, and Wiegele \cite{MPRW} further showed that
\reff{sdp:augLag:1st} can be reduced to \be  \label{sdp:augLag}
\baray{rl} \underset{y\in\re^m}{\rm max} &
 b^Ty - \frac{\sigma}{2} \|(\mc{A}^*(y)-C+Y/\sigma)_{\mc{K}}\|^2
+ \frac{1}{2\sigma} \|Y\|^2. \earay \ee

When $\mc{K}=\mc{S}_+^N$ is a single product, Malick, Povh, Rendl,
and Wiegele \cite{MPRW} proposed a general framework
(cf. \cite[Algorithm~4.3]{MPRW}) of regularization methods for solving
large scale SDP problems. It can be readily generalized to the case
that $\mc{K}$ is a product of several semidefinite cones. We
describe it as follows:
\begin{alg}
\label{regalg:MPRW}
Choose $\eps^{in},\eps^{out} \in (0,1)$.  \\
\indent Initialization: Choose $Y_0 \in \mc{S}^N$, $Z_0= 0$, $y_0\in \re^m$, $\sig_0$ and set $k=0$. \\
\indent While ($\|Z_k+\mc{A}^*(y_k)-C\| \geq \eps^{out}$) ({\it outer loop}): \\
\indent  \indent  Set $j:=0$, $y_{k,j} := y_k$ and $X_{k,j} := Y_k$. \\
\indent \indent  While ($\|b-\mc{A}(X_{k,j})\| \geq \eps^{in}$) ({\it inner loop}): \\
\indent \indent \indent  Compute the projections \\
\indent \indent \indent \indent  $X_{k,j} :=
\sigma_k(Y_k/\sigma_k+\mc{A}^*(y_{k,j})-C)_{ \mc{K} }$,\\
\indent \indent \indent \indent  $Z_{k,j} :=
-(Y_k/\sigma_k+\mc{A}^*(y_{k,j})-C)_{-\mc{K}}. $\\
\indent \indent \indent Set $g_j := b-\mc{A}(X_{k,j})$.\\
\indent \indent \indent Set $y_{k,j+1} :=  y_{k,j} + \tau W g_j$ with appropriate $\tau$ and $W$. \\
\indent \indent  \indent  Set $j:=j+1$. \\
\indent \indent  end ({\large \it inner loop}) \\
\indent \indent Set $Y_{k+1}:=X_{k,j}$, $y_{k+1} := y_{k,j}$ and update $\sigma_{k}$. \\
\indent  \indent  Set $k:=k+1$. \\
\indent end ({\large \it outer loop})
\end{alg}

In Algorithm~\ref{regalg:MPRW}, if $W$ is chosen to be
$(\mc{A}\mc{A}^*)^{-1}$ and $\tau = 1/\sigma$,
Algorithm~\ref{regalg:MPRW} becomes the {\it Boundary Point
Method (BPM)}, which was originally proposed by Povh, Rendl, and
Wiegele \cite{PRW06} (also see \cite{MPRW}) for solving big SDP problems.
This method was proven efficient in some applications, as
illustrated in \cite{MPRW,PRW06}. The description of this method is:
\begin{alg} \label{BPM:Ksdp}
Choose $\eps \in (0,1)$.  \\
\indent Initialization: Choose $Y_0 \in \mc{S}^N$, $Z_0= 0$,
$y_0\in \re^m$, $\sig_0$ and set $k:=0$. \\
\indent While ($\|b-\mc{A}(X_k)\| \geq \eps$ or $\|C-Z_k-\mc{A}^*(y_k)\| \geq \eps$) \\
\indent \indent Solve $\mc{A}\mc{A}^*y_{k+1} = \mc{A}(C-Z_k) + (b-\mc{A}(Y_k))/\sigma_{k}$
for $y_{k+1}$.\\
\indent \indent Compute the projections \\
\indent \indent \indent  $X_{k+1} :=
\sigma_k(Y_k/\sigma_k+\mc{A}^*(y_{k+1})-C)_{\mc{K}}$, \\
\indent \indent \indent   $Z_{k+1} :=
-(Y_k/\sigma_k+\mc{A}^*(y_{k+1})-C)_{-\mc{K}}. $\\
\indent \indent Set $Y_{k+1}:=X_{k+1}$ and update $\sigma_{k}$.  \\
\indent \indent Set $k: = k+1.$ \\
 \indent end
\end{alg}

For solving large scale SDP relaxations in polynomial optimization,
Algorithm~\ref{BPM:Ksdp} might converge fast at the beginning, but
generally has relatively slow convergence when it gets close to
optimal solutions. This is because it is basically a gradient type
method. When eigenvalue decompositions are not expensive,
Algorithm~\ref{BPM:Ksdp} usually works very well, as demonstrated in
\cite{MPRW,PRW06}. When it is expensive to compute eigenvalue
decompositions, if Algorithm~\ref{BPM:Ksdp} takes a large number of
iterations, then it might consume a lot of time. On the other hand,
Algorithm~\ref{BPM:Ksdp} has simple iterations and is easily
implementable. In each iteration, we only need to solve a linear
system (its coefficient matrix $\mc{A}\mc{A}^*$ is fixed) and
compute an eigenvalue decomposition. This is a big advantage. It
would be applied to get a good approximate solution.

To get more accurate solutions, we need more efficient methods for
the inner loop of Algorithm~\ref{regalg:MPRW}. Typically, Newton
type methods have good properties like local superlinear or
quadratic convergence. This leads to the {\it Newton-CG
Augmented Lagrangian} method, which was proposed by Zhao, Sun and
Toh \cite{ZST08}. It also has good numerical
performance in solving big SDP problems. In the following, we
describe this important method for solving
\reff{Ksdp:prm}-\reff{Ksdp:dual} when $\mc{K}$ is a product of
several semidefinite cones.

Denote by $\varphi_{\sigma}(Y,y)$ the objective in \reff{sdp:augLag}.
When $\mc{K}=\mc{S}_+^N$ is a single product, $\varphi_{\sigma}(Y,y)$ is differentiable
\cite[Proposition~3.2]{MPRW} and
\[
\nabla_y \varphi_{\sigma}(Y,y) = b - \sigma
\mc{A}\Big((\mc{A}^*(y)-C+Y/\sigma)_{\mc{K}}\Big).
\]
The above is also true when $\mc{K}$ is a product of several
semidefinite cones. The inner loop of Algorithm~\ref{regalg:MPRW} is
solving the maximization problem \be  \label{max:varphi}
\max_{y\in\re^m}  \quad \varphi_{\sigma}(Y_k,y). \ee Since
$\varphi_{\sigma}$ is concave in $y$, a point $\hat{y}$ is a
maximizer of \reff{max:varphi} if and only if
\[
\nabla_y \varphi_{\sigma}(Y_k,\hat{y})=0.
\]
The function $\varphi_{\sigma}(Y,y)$ is not twice differentiable, so
the standard Newton's method is not applicable. However, the
function $\varphi_{\sigma}(Y,y)$ is semismooth, and semismooth
Newton's method could be applied to get local superlinear or
quadratic convergence, as pointed out in \cite[Section~3.2]{ZST08}.
For this purpose, we need the generalized Hessian of $\varphi_\sig$ in computation.
We refer to \cite[Section~3.2]{ZST08} for a numerical method of
evaluating $\nabla_y^2 \varphi_{\sigma}(Y,y)$.
It is important to point out that the Hessian $\nabla_y^2
\varphi_{\sigma}(Y,y)$ does not need to be explicitly formulated. It
is implicitly available such that the matrix vector product
$\nabla_y^2 \varphi_{\sigma}(Y,y) \cdot z$ can be evaluated efficiently.
Generally, we always have $\nabla_y^2 \varphi_{\sigma}(Y,y) \succeq
0$, and $\nabla_y^2 \varphi_{\sigma}(Y,y) \succ 0$ if some
nondegeneracy condition holds (cf. \cite[Prop.~3.2]{ZST08}). In either
case, an approximate semismooth Newton direction $d_{new}$ for
\reff{max:varphi} can be determined from the linear system \be
\label{linsys:Newton} \Big( \nabla_y^2 \varphi_{\sigma}(Y,y) + \eps
\cdot I_N \Big)  d_{new} = \nabla_y  \varphi_\sigma. \ee Here
$\eps>0$ is a tiny number ensuring the positive definiteness of the
above linear system. When $m$ is huge, it is usually not practical
to solve \reff{linsys:Newton} by direct methods like Cholesky
factorization. To avoid this difficulty, conjugate gradient (CG)
iterations are suitable, as proposed in \cite{ZST08}.
%

Now we describe the Newton-CG Augmented Lagrangian regularization method.
\begin{alg} \label{alg:Newton-CG}
Choose $\eps^{in},\eps^{out} \in (0,1)$, $\eps>0$, $\dt\in(0,1)$, $\rho>1$, $\sigma_{\max},K\in \N$. \\
\indent Initialization: Choose $X_0,Z_0\in \mc{S}^N$,
$y_0\in \re^m$, $\sig_0$ and set $k:=0$. \\
\indent While ($\|Z_k+\mc{A}^*(y_k)-C\| \geq \eps^{out}$) ({\large \it outer loop}): \\
\indent  \indent  Set $Y_k:=X_k$. \\
\indent  \indent  Set $j:=0$ and $y_{k,j} := y_k$. \\
\indent \indent While ($\|\nabla_y \varphi_{\sigma_k}(Y_k,y_{k,j})\|
\geq \eps^{in}$)
({\large \it inner loop}): \\
\indent  \indent \indent Set $g_j := \nabla_y \varphi_{\sigma_k}(Y_k,y_{k,j})$. \\
\indent  \indent \indent Compute $d_{new}$ from \reff{linsys:Newton}
by applying preconditioned CG at most $K$ steps. \\
\indent  \indent \indent  Find the smallest integer $\af>0$ such that
\be \label{cond:linsch}
\varphi_{\sigma_k}(Y_k, y_{k,j}+ \dt^\af \cdot
d_{new} ) \geq \varphi_{\sigma_k}(Y_k, y_{k,j}) + \dt^\af \cdot
g_j^T d_{new}. \ee
\indent  \indent \indent  Set $y_{k,j+1} := y_{k,j} + \dt^\af \cdot d_{new}$. \\
\indent \indent  \indent  Compute the projections: \\
\indent \indent  \indent  \indent $X_k :=
\sigma_k(Y_k/\sigma_k+\mc{A}^*(y_{k,j+1})-C)_{\mc{K}}$,\\
\indent \indent  \indent  \indent $Z_k := -(Y_k/\sigma_k+\mc{A}^*(y_{k,j+1})-C)_{-\mc{K}}$.\\
\indent \indent  \indent  Set $j:=j+1$. \\
\indent \indent  end ({\large \it inner loop}) \\
\indent  \indent  Set $y_{k+1} := y_{k,j}$. \\
\indent  \indent  If $\sigma_k \leq \sigma_{\max}$, set $\sigma_{k+1} := \rho \sigma_k$. \\
\indent  \indent  Set $k:=k+1$. \\
\indent end ({\large \it outer loop})
\end{alg}

When $\mc{K}=\mc{S}_+^N$ is a single product, the convergence of
Algorithm~\ref{alg:Newton-CG} has been discussed in
\cite[Theorems~3.5, 4.1, 4.2]{ZST08}. These results could be readily
generalized to $\mc{K}$ being a product of several semidefinite cones.
The specifics about the convergence are beyond the scope of this
paper. We refer to \cite{MPRW,Roc76CO,Roc76MOR,ZST08}.

\section{Computational experiments}
\setcounter{equation}{0}

This section presents numerical experiments of applying
Algorithm~\ref{alg:Newton-CG} in solving polynomial optimization problems. An
excellent implementation of Algorithm~\ref{alg:Newton-CG} is
software {\tt SDPNAL} \cite{sdpnal}. We use it to solve the SDP
relaxations (its earlier version posted in early 2010 was used).
The computation is implemented in Matlab 7.10
on a Dell Linux Desktop running CentOS (5.6) with 8GB
memory and Intel Core CPU 2.8GHz.
We use the default parameters of {\tt
SDPNAL}: $\sigma_0 =10,~K=500,~ {\tt Tol} = 10^{-6}$. Set
\be \label{def:rp:rd}
R_P :=\frac{\|\mathcal{A}(X)-b\|_2}{1+\|b\|_2},~~
R_D := \frac{\|\mathcal{A}^{*}(y)+Z-C\|_2}{1+\|C\|_2}, \nonumber
\ee
which measure the feasibilities of the computed solutions for the
primal and dual SDP problems respectively. We terminate computation
when $\max\{R_P,R_D\}\leq {\tt Tol}$. Other parameters are set to be
the default ones of {\tt SDPNAL}.

If the computed dual optimal solution $Z^*$ of \reff{sdp:momsos} or
\reff{dualsdp:consos} satisfies FEC, we could
extract a global minimizer $x^*$;
otherwise, we just set $Z^*(2:n+1,1)$ as a starting
point and get a local optimal solution $x^{*}$ by using nonlinear
programming solvers in Matlab Optimization Toolbox. In
either case, the error of computed $x^*$ is measured as
\be \label{df:sol-err}
\mbox{errsol} \quad = \quad \frac{|f(x^*)-\underline{f}|}{\max\{1,
|f(x^*)|\} }, \ee
where $\underline{f}$ is a lower bound returned by solving
the SDP relaxation. The error of a computed optimal triple $(X,y,Z)$
for the SDP relaxation itself is measured as
\be \label{def:sdp-err}
\mbox{errsdp} \quad = \quad \max\left\{ \frac{|b^{\top}y-\langle
C,X\rangle|}{1+|b^{\top}y|+|\langle C,X\rangle|}, R_P, R_D\right\}.
\ee
The consumed computer time is in the format {\tt hr:mn:sc} with hr
(resp. {\tt mn, sc}) standing for the consumed hours (resp. minutes,
seconds). In the tables of this paper, {\tt min, med} and {\tt max}
respectively stands for the minimum, median, and maximum of
quantities like time, errsol, errsdp, etc.

We would like to point out that
extracting global minimizers is a difficult problem.
When FEC is satisfied, one or several global minimizers of \reff{pop:S}
can be found by solving eigenvalue problems (cf. \cite{HL05}).
When FEC fails, it's an open question how to extract global
minimizers. In such situations, we just use nonlinear programming
methods to get a local minimizer with $Z^*(2:n+1,1)$ being a starting
point. The experiment results in \cite{WKKM} show that this approach
usually works very well. In our experiments, we also use this
technique to get a local minimizer when FEC fails.


The testing problems in our
experiments are in three categories: (a) unconstrained polynomial
optimization and it's application in sensor network localization;
(b) homogeneous polynomial optimization and it's application in
computing stability numbers; (c) constrained polynomial
optimization.

\subsection{Unconstrained polynomial optimization}
\label{subsec:uc-num}

\begin{exm}
Minimize the following least squares polynomial
\[
\sum_{k=1}^3  \left(\sum_{i=1}^n x_i^k - 1\right)^2  + \sum_{i=1}^n
\left(x_{i-1}^2+x_i^2+x_{i+1}^2 - x_i^3-1\right)^2
\]
where $x_0=x_{n+1}=0$. For $n=16$, the resulting SDP relaxation
\reff{sdp:stdsos}-\reff{sdp:momsos} has size $(N,m)=(969, 74612)$.
Solving it by {\tt SDPNAL} takes about 34 minutes. The computed
solution of the SDP relaxation has error around $6\cdot 10^{-7}$.
The computed lower bound $f_{sos}^{uc} \approx 7.5586$. The optimal
$Z^*$ has rank two and FEC holds. We get two optimal solutions.
Their errors are about $2\cdot 10^{-6}$. \qed
\end{exm}

\begin{table}[htb]
\centering
\begin{scriptsize}
\begin{tabular}{|r||c||c||r|r|r||c||c|} \hline
$n$ & (N,m) & \#Inst   & \multicolumn{3}{c||}{time (min, med, max)}
& errsol (min, med, max)& errsdp (min, med, max) \\ \hline 20 &
(231, 10625) & 20 &   0:00:02 &
0:00:04   &   0:00:09  & (4.1e-7, 1.0e-5, 1.6e-4)&(2.5e-7, 7.3e-7, 1.3e-6) \\
\hline 30 & (496,
46375) & 20 & 0:00:12 & 0:00:21  & 0:00:31 & (1.3e-7, 5.6e-5, 1.5e-4)& (3.2e-7, 6.8e-7, 1.0e-6) \\
\hline 40 &  (861, 135750) & 10 & 0:00:57 & 0:01:10 & 0:01:24 &
(7.8e-7, 1.2e-4, 3.1e-4)& (4.2e-7, 4.7e-7, 9.6e-7)\\ \hline 50 &
(1326, 316250) & 5 &
0:02:44 & 0:03:17 & 0:04:08 & (1.3e-5, 3.2e-5, 2.3e-4)& (5.6e-7, 6.4e-7, 8.3e-7) \\
\hline 60 & (1891,
635375) & 5 & 0:07:55&  0:08:49& 0:09:48& (4.6e-5, 1.8e-4, 5.1e-4)& (4.8e-7, 9.1e-7, 9.5e-7)\\
\hline 70 & (2556, 1150625) & 5 & 0:17:38 & 0:19:34 & 0:22:33 &
(8.0e-5, 2.8e-4, 3.3e-4)&(4.1e-7, 5.7e-7, 9.2e-7) \\ \hline 80 &
(3321, 1929500) & 3 & 0:38:45 & 0:38:48 & 0:42:46  & (9.3e-5,
2.7e-4, 9.6e-4) & (3.7e-7, 7.0e-7, 9.9e-7)
\\ \hline 90 & (4186,
3049500) & 3 & 1:37:04 & 1:46:57 & 2:02:01 & (1.1e-4, 2.7e-4, 6.4e-4) & (4.3e-7, 5.2e-7, 9.5e-7)\\
\hline 100& (5151, 4598125) & 3 & 2:48:03 & 2:55:34 & 3:35:27 &
(2.1e-4, 2.6e-4, 4.5e-4) &  (7.1e-7, 7.9e-7, 8.7e-7)\\ \hline
\end{tabular}
\end{scriptsize}
\caption{Computational results for random unconstrained optimization of
degree $4$} \label{UCtab:rand:deg=4}
\end{table}

\begin{table}[htb]
\centering
\begin{scriptsize}
\begin{tabular}{|r||c||c||r|r|r||c||c|} \hline
$n$ & (N,m) & \#Inst   & \multicolumn{3}{c||}{time (min, med, max)}
& errsol (min, med, max)& errsdp (min, med, max)  \\ \hline 10 &
(286, 8007) & 20 & 0:00:07  &
0:00:17  & 0:00:36     & (2.7e-7, 5.2e-6, 6.6e-5)&(2.4e-8, 4.2e-7, 1.1e-6) \\
\hline 15 & (816, 54263) & 10  &  0:01:12 &  0:01:51    & 3:07:37 &
 (5.1e-6, 3.6e-5, 7.0e-5) & (2.0e-7, 7.4e-7, 9.6e-7) \\ \hline
20 &  (1771, 230229)  & 3  &  2:54:42    &  2:57:57    & 15:10:08 &
  (1.4e-4, 2.5e-4, 4.0e-4) & (3.1e-7, 4.9e-7, 6.0e-7) \\ \hline
25 &  (3276, 736280)  & 3  &   2:02:59  &  5:25:06  & 7:34:03 &
(1.6e-3,  8.0e-3, 4.7e-2) & (2.6e-6,  8.6e-6, 5.7e-5)
 \\ \hline
\end{tabular}
\end{scriptsize} \caption{Computational results for random
unconstrained optimization of degree $6$ } \label{UCtab:rand:deg=6}
\end{table}

\begin{table}[htb]
\centering \begin{scriptsize}
\begin{tabular}{|r||c||c||r|r|r||c||c|} \hline
$n$ & (N,m) & \#Inst  & \multicolumn{3}{c||}{time (min, med, max)} &
errsol (min, med, max)& errsdp (min, med, max) \\ \hline 8 &   (495,
12869) & 20  &  0:00:18& 0:00:42&  0:01:11   &
  (1.6e-7, 2.6e-5, 5.6e-4) & (1.0e-7, 5.8e-7, 4.1e-6)\\ \hline
10 &  (1001, 43757)  & 10  & 0:04:46 & 0:06:42 &  0:08:05   &
(3.9e-5, 8.4e-5, 5.3e-4)& (2.4e-7, 5.9e-7, 3.0e-6)\\ \hline 12 &
(1820, 125969)  & 3 & 0:26:32 &  0:37:43 & 1:02:37 & (1.3e-5, 4.4e-5, 5.7e-3)&(1.1e-7, 7.2e-7, 5.3e-6) \\
\hline 15 &  (3876, 490313) & 3 & 6:31:11     & 9:08:29
& 10:21:21 & (6.8e-4, 8.1e-4, 4.5e-3) &(9.9e-7, 1.1e-6, 5.6e-6) \\
\hline
\end{tabular}
\end{scriptsize}
\caption{Computational results for random unconstrained optimization of
degree $8$ } \label{UCtab:rand:deg=8}
\end{table}

\begin{table}[htb]
\centering
\begin{scriptsize}
\begin{tabular}{|r||c||c||r|r|r||c||c|} \hline
$n$ & (N,m) & \#Inst  & \multicolumn{3}{c||}{time (min, med, max)} &
errsol (min, med, max) & errsdp (min, med, max)\\ \hline 6 &  (462,
8007) & 20 &  0:00:10 & 0:00:18 & 0:00:32    &
  (3.6e-7, 1.4e-5, 1.4e-4)&(3.2e-8, 5.2e-7, 3.1e-6) \\ \hline
8 &  (1287, 43757)  & 10  &  0:04:13 & 0:05:33    & 0:10:23 &
  (5.6e-6, 5.2e-5, 3.1e-4)& (2.2e-7, 3.9e-7, 1.8e-6) \\ \hline
9 &  (2002, 92377)  & 3 & 0:13:13 &  0:18:31  & 0:43:28 &
 (2.2e-4, 7.3e-4, 8.4e-4) & (1.1e-6, 2.5e-6, 2.9e-6)\\ \hline
10 & (3003, 184755)  & 3  &   3:53:13 &   3:58:15  & 4:02:11 &
  (2.3e-3, 2.4e-3, 4.1e-3)&(4.7e-7, 1.2e-6, 4.2e-6) \\ \hline
\end{tabular}
\end{scriptsize} \caption{Computational results for random
unconstrained optimization of degree $10$ } \label{UCtab:rand:deg=10}
\end{table}

\begin{exm}[Random polynomials]
We test the performance of Algorithm~\ref{alg:Newton-CG}
(via {\tt SDPNAL}) in solving SDP relaxations for randomly generated
polynomial optimization problems. To ensure the existence of a global
minimizer, generate $f(x)$ randomly as
\[
f(x) = f^T[x]_{2d-1} + [x^d]^T F^TF [x^d],
\]
where $f/F$ is a Gaussian random vector/matrix of a proper
dimension. Here $[x^d]$ denotes the vector of monomials of degree
equal to $d$.
The computational results are shown in
Tables~\ref{UCtab:rand:deg=4}-\ref{UCtab:rand:deg=10}. There
$\#$Inst denotes the number of randomly generated instances, and
$(N,m)$ denotes the size of the corresponding SDP relaxation
\reff{sdp:stdsos}-\reff{sdp:momsos}.

When $f(x)$ has degree $4$ ($d=2$), SDP relaxation
\reff{sdp:stdsos}-\reff{sdp:momsos} is solved quite well. For
$n=20\sim 30$, the computation takes up to half a minute; for
$n=40\sim 60$, it takes a couple of minutes; for $n=70\sim 80$, it
takes less than one hour; for $n=90\sim 100$, it takes a few hours.
When $f(x)$ has degree $6$ ($d=3$), for $n=15$, solving
\reff{sdp:stdsos}-\reff{sdp:momsos} takes up to a few hours; for
$n=20\sim 25$, it takes a couple of hours.
When $f(x)$ has degree $8$ ($d=4$), for $n=10$, solving
\reff{sdp:stdsos}-\reff{sdp:momsos} takes a couple of minutes; for
$n=12\sim 15$, it takes about one to ten hours.
When $f(x)$ has degree $10$ ($d=5$),  for $n=8$, solving
\reff{sdp:stdsos}-\reff{sdp:momsos} takes a couple of minutes; for
$n=9$, it takes less than one hour; for $n=10$, it takes a few
hours.
From Tables~\ref{UCtab:rand:deg=4} to \ref{UCtab:rand:deg=10}, we
can see that the SDP relaxations are solved successfully. The
obtained solutions for polynomial optimization are also reasonably
very well. They are slightly less accurate than the computed
solutions of the SDP relaxation itself. This is probably
because the SDP relaxation \reff{sdp:stdsos}-\reff{sdp:momsos} is
not exact in minimizing the generated polynomials.

%
%

The computations here show that Algorithm~\ref{alg:Newton-CG} could
solve large scale polynomial optimization problems. A quartic polynomial
optimization with $100$ variables could be solved within a couple of
hours on a regular computer. This is almost impossible by using SDP
solvers based on interior point methods. \qed \end{exm}

\begin{exm}[Sensor Network Localization]
Given a graph $G=(V,E)$ and a distance for each edge, the
sensor network localization problem
is to find locations of vertices so that their distances are equal
to the desired ones. This problem can be formulated as follows: find
a sequence of unknown vectors ({\it sensors}) $u_1,u_2,\ldots,u_s
\in \re^k$ (typically $k=1,2,3$, we focus on $k=2$ in this example)
such that the distances between these sensors and some other fixed
vectors ({\it anchors}) $a_1,\ldots,a_\ell$ are equal to given
distances. Recently, there is much work on solving
sensor network localization by SDP
techniques, like \cite{BisYe,KKW09,Nie09}. Given edge subsets
\[
\mc{E}_{S}\subset \{(i,j):\, 1 \leq i<j \leq s \}, \quad
\mc{E}_{A}=\{(i,j): 1\leq i \leq s, 1\leq j \leq \ell \},
\]
for every $(i,j) \in \mc{E}_{S}$, let $d_{ij}$ be the distance
between $u_i$ and $u_j$, and for every $(i,j) \in \mc{E}_{A}$, let
$e_{ij}$ be the distance between $u_i$ and $a_j$. Denote $u_i =
(x_{ki-k+1},\ldots, x_{ki})$ for $i=1,\ldots,s$. The
sensor network localization problem
is equivalent to finding coordinates $x_{k1},\ldots,x_{ks}$
satisfying the equations
\begin{align*}
\|u_i - u_j \|_2^2 
=  d_{ij}^2\,\quad \forall\, (i,j) \in\mc{E}_{S}, \qquad
\|u_i - a_j \|_2^2 
= e_{ij}^2\, \quad \forall\, (i,j) \in\mc{E}_{A}.
\end{align*}
It is also equivalent to the quartic polynomial optimization problem
\begin{align} \label{L2err}
\min_{u_1,\ldots, u_s} \,\quad\, \sum_{(i,j)\in \mc{E}_{S}} \left(
\| u_i - u_j \|_2^2 - d_{ij}^2\right)^2 +\, \sum_{(i,j)\in
\mc{E}_{A}} \left( \| u_i - a_j \|_2^2 - e_{ij}^2 \right)^2.
\end{align}
Typically, it is large scale. We use {\tt SDPNAL} to solve its SDP
relaxation \reff{sdp:stdsos}-\reff{sdp:momsos}.
\begin{table}[htb]
\centering
\begin{scriptsize}
\begin{tabular}{|c||c||r|r|r||c||c|} \hline
\#sensor   & \#Inst   & \multicolumn{3}{c||}{time
(min, med, max)} & RMSD (min, med, max)& errsdp (min, med, max) \\
\hline 15  & 15 & 0:00:24 & 0:00:52   &   0:02:02  & (8.1e-6,
2.4e-5, 1.4e-4)&(1.1e-7, 4.2e-7, 1.6e-6) \\\hline 20  & 15 & 0:02:04
& 0:03:19   &   0:09:12  & (1.5e-5, 5.5e-5, 1.5e-4)&(2.9e-7, 4.4e-7,
2.0e-6) \\\hline 25  & 10 & 0:14:18 &
0:35:02   &   1:12:21  & (4.3e-5, 8.7e-5, 2.2e-4)&(2.4e-7, 6.3e-7, 1.6e-6) \\
\hline 30   &10 & 1:22:18 & 2:44:05   &   5:51:36  & (2.3e-5,
2.3e-4, 2.7e-3)&(9.2e-8, 1.8e-6, 5.3e-4) \\ \hline 35  &  3 &
09:59:35 & 19:13:58   &  27:08:37 & (1.3e-3, 1.6e-3,
2.2e-3)&(6.5e-6, 5.1e-5, 6.5e-4) \\ \hline 40  &  3 & 48:33:59 &
50:54:34   &
61:19:58 & (1.2e-3, 1.6e-3, 2.7e-3)&(2.2e-3, 3.2e-3, 4.0e-3) \\
\hline
\end{tabular}
\end{scriptsize}
\caption{Computational results for sensor network localization problems.}
\label{sensor:network:exact}
\end{table}
To test its performance, we randomly generate sensors $u_1,\ldots,
u_s$ from the square $[-0.5,\,0.5]\times [-0.5,\,0.5]$. Fix four
anchors as $(\pm 0.45,\, \pm 0.45)$. For each pair $(i,j)$, select
it to $\mathcal{\mc{E}_{S}}$ with probability $0.6$ and to
$\mathcal{\mc{E}_{A}}$ with probability $0.3$. Then compute each
distance $d_{ij}$ and $e_{ij}$. After the SDP relaxation is solved,
we use $Z^*(2:n+1,1)$ as a starting point and apply function
{\tt lsqnonlin} in Matlab Optimization Toolbox to get a local minimizer
$(\hat{u}_1,\ldots, \hat{u}_s)$ of \reff{L2err} (we use the
technique that was proposed in \cite{KKW09}). The errors of computed locations are
measured by the Root Mean Square Distance $ \text{RMSD}
=(\frac{1}{s}\Sig_{i=1}^{s}\|\hat{u}_i-u^{*}_{i}\|^2)^{1/2}, $ as
used in \cite{BisYe}.

The computational results are shown in Table~\ref{sensor:network:exact}.
We can see that the SDP
relaxation of \reff{L2err} is solved reasonably well. In many
instances, FEC is not satisfied,
so we can only get a local minimizer of \reff{L2err}
by using the technique from \cite{KKW09}.
The true locations of sensors are found
with small errors. Possible reasons for FEC fails might be:
the SDP relaxation was not solved accurately enough,
or it is not exact for \reff{L2err}.

We would like to remark that the SDP
relaxation \reff{sdp:stdsos}-\reff{sdp:momsos} for \reff{L2err}
does not exploit the sparsity pattern.
There exists work of using sparse SDP or SOS type relaxations
for solving sensor network localization problems
(e.g., Kim et al. \cite{KKW09} and Nie \cite{Nie09}).
Generally, solving \reff{sdp:stdsos}-\reff{sdp:momsos} for \reff{L2err}
is much more difficult than solving its sparse versions like in \cite{KKW09,Nie09}.
The numerical experiments in \cite{KKW09,Nie09} show that
exploiting sparsity will allow us to solve much bigger problems.
However, in the view of quality of approximations,
sparse SDP relaxations are typically weaker than the general dense one.
Thus, it is still meaningful if we can solve \reff{sdp:stdsos}-\reff{sdp:momsos} for
large scale sensor network localization problems.
\qed
\end{exm}

\subsection{Homogeneous polynomial optimization}

\begin{exm} Minimize the following square free quartic form over $\mathbb{S}^{n-1}$
\[
\sum_{ 1 \leq i < j < k < \ell \leq n}  (-i-j+k+\ell)x_ix_jx_kx_\ell
.
\]
For $n=50$, the resulting SDP \reff{hmgsdp:prime}-\reff{hmgsdp:dual}
has size $(N,m)=(1275,292824)$. Solving
\reff{hmgsdp:prime}-\reff{hmgsdp:dual} takes about $38$ minutes. The
error of the computed solution for the SDP relaxation is around
$8\cdot 10^{-8}$. The computed lower bound $f_{sos}^{hmg} \approx
-140.4051$. The optimal $Z^*$ has rank two and FEC holds, so we get
two optimizers.
Their errors are around $2\cdot 10^{-7}$.
\qed \end{exm}

\begin{exm} Minimize the following sextic form over $\mathbb{S}^{n-1}$
\[
\sum_{1\leq i \leq n} x_i^6 + \sum_{1\leq i \leq n-1}x_i^3x_{i+1}^3.
\]
For $n=20$, the size of the resulting SDP
\reff{hmgsdp:prime}-\reff{hmgsdp:dual} is $(N,m)=(1540,177099)$. In
this problem, we set parameter ${\tt
Tol}=10^{-10}$ in running {\tt SDPNAL} (For the default choice ${\tt
Tol}=10^{-6}$, {\tt SDPNAL} does not converge very well for this
example), i.e., we terminate the
computation when $\max\{R_p,R_D\}<10^{-10}$. Solving
\reff{hmgsdp:prime}-\reff{hmgsdp:dual} takes about $2.5$ hours. The
computed solution of the SDP relaxation has error around $1\cdot
10^{-5}$. The computed lower bound $f_{sos}^{hmg} \approx
1.1451\times 10^{-4}$.
The computed optimal $Z^*$ has rank one and we get one global
minimizer from it.
Its error is around $1.3 \cdot 10^{-5}$.
\qed \end{exm}

\begin{exm} Minimize the following sextic form over $\mathbb{S}^{n-1}$
\[
\sum_{ 1 \leq i < j < k \leq n} x_i^2x_j^2x_k^2 + x_i^3x_j^2x_k +
x_i^2x_j^3x_k + x_ix_j^3x_k^2.
\]
For $n=20$, the resulting SDP \reff{hmgsdp:prime}-\reff{hmgsdp:dual}
has size $(N,m)=(1540,177099)$. Solving
\reff{hmgsdp:prime}-\reff{hmgsdp:dual} takes about $1.8$ hours. The
error of the computed solution of the SDP relaxation is around
$1.7\cdot 10^{-6}$. The computed lower bound $f_{sos}^{hmg} \approx
-0.3827$. The computed optimal $Z^*$ has rank one, so we get one
global minimizer. 
%
%
Its error is around $7.4\cdot 10^{-7}$. \qed
\end{exm}

An important application of homogenous polynomial optimization
\reff{hmgpop:sphere} is computing stability numbers of graphs.

\begin{exm}[Stability numbers of graphs]
Let $G=(V,E)$ be a graph with $|V|=n$. The stability number $\af(G)$
is the cardinality of the biggest stable subset(s) (their vertices
are not connected by any edges) of $V$. It was shown in Motzkin and
Straus \cite{MS65} (also see De Klerk and Pasechnik \cite{dKP02})
that
\[
\af(G)^{-1} = \min_{x\in \Dt_n} \quad x^T(A+I_n)x,
\]
where $\Dt_n$ is the standard simplex in $\re^n$ and $A$ is the
adjacency matrix associated with $G$. If replacing every $x_i\geq 0$
by $x_i^2$, we get \be \label{af(G):M-S} \af(G)^{-1} = \min_{\|x\|_2
= 1} \quad \sum_{i=1}^n x_i^4 + 2 \sum_{(i,j)\in E} x_i^2x_j^2 . \ee
This is a quartic homogeneous polynomial optimization. When a lower
bound $f_{sos}^{hmg}$ of \reff{af(G):M-S} is computed from its SDP relaxation,
we round $\left(f_{sos}^{hmg}\right)^{-1}$ to the nearest
integer which will be used to estimate $\af(G)$.

We generate random graphs $G$, and solve the SDP relaxation of
\reff{af(G):M-S}. The generation of random graphs is in a similar
way as in Bomze and De Klerk \cite[Section~6]{BdK02}. For $n=20, 30,
40, 50, 60$, we generate random graphs $G=(V,E)$ with $|V|=n$.
Select a random subset $M\subset V$ with $|M|=n/2$. The edges
$e_{ij}(\{i,j\} \not\subset M)$ are generated with probability
$\half$. The computational results are in Table \ref{stab:test}.
\begin{table}
\centering
\begin{scriptsize}
\begin{tabular}{|r||c||c||r|r|r||c|} \hline
$n$ & (N,m) &\# Inst   & \multicolumn{3}{c||}{time (min, med, max)}& errsdp (min, med, max) \\
\hline 20 &  (210, 8854)  & 20  & 0:00:06 & 0:00:11  & 0:00:26 & (1.2e-7, 6.7e-7, 1.0e-6)\\
\hline 30 &  (465, 40919) & 20  & 0:00:45 & 0:01:21  & 0:01:52 & (2.0e-7, 5.2e-7, 1.1e-6) \\
\hline 40 &  (820, 123409)& 10  & 0:02:31 & 0:05:19  & 0:08:58 & (4.1e-7, 7.4e-7, 1.6e-6)\\
\hline 50 & (1275, 292824)& 10  & 0:13:30 & 0:19:10 & 0:29:29  & (4.5e-7, 5.9e-7, 9.6e-7)\\
\hline 60 & (1830, 595664)& 5  & 0:44:19 & 1:05:32 & 1:47:51 & (2.4e-7, 5.6e-7, 3.8e-6) \\
\hline 70 & (2485, 1088429)& 5 &  2:33:24 &  4:20:13  &  5:07:37 & (4.2e-7, 6.2e-7, 7.9e-7)\\
\hline 80 & (3240, 1837619)& 3  & 7:31:21 & 9:43:40 & 10:52:27 &(3.6e-7, 4.4e-7, 7.8e-7) \\
\hline 90 & (4095, 2919734)& 3  & 17:10:41 & 17:44:02 & 18:45:28 &(2.1e-7, 3.7e-7, 3.7e-7)\\
\hline
\end{tabular}
\end{scriptsize} \caption{Computational results for stability
number of random graphs } \label{stab:test}
\end{table}
As one can see, for $n=20,30,40,50$, solving
\reff{hmgsdp:prime}-\reff{hmgsdp:dual} takes less than half an hour;
for $n=60,70$, it takes a few hours; for $n=80,90$, it takes $7$ to
$19$ hours. In all the instances, we get correct stability numbers.
All the SDP relaxations themselves are also solved successively.
\qed
\end{exm}

\subsection{Constrained polynomial optimization}

\begin{exm}
Minimize the sextic polynomial
\begin{align*}
\sum_{i=1}^n ix_i^3 + \sum_{1 \leq i < j \leq n} (i+j) x_i^3x_j^3
\end{align*}
over the unit ball $B(0,1)$. We apply the $3^{rd}$ Lasserre's
relaxation \reff{consos:deg2d}. The resulting cone $\mc{K}$ has $2$
blocks. When $n=20$, solving
\reff{prmsdp:consos}-\reff{dualsdp:consos} takes about $34$ minutes.
The computed lower bound $f_{sos}^{con}=-20$. The optimal $Z^*$ has
rank one, and we get the solution:
\[
(0,   0,   0,   0,   0,   0,   0,   0,   0,   0,   0,   0,
 0,   0,   0,   0,   0,   0,   0,  -1).
\]
It is feasible and a global minimizer.
\qed \end{exm}

\begin{exm}
Consider the polynomial optimization
\begin{align*}
\min_{x\in\re^n}  & \quad \sum_{1 \leq i < j <k \leq n/2} \left(
x_ix_jx_k+ x_{\frac{n}{2}+i}x_{\frac{n}{2}+j}x_{\frac{n}{2}+k} +
x_ix_jx_kx_{\frac{n}{2}+i}x_{\frac{n}{2}+j}x_{\frac{n}{2}+k} \right) \\
s.t. & \quad x_1^4+\cdots+x_{\frac{n}{2}}^4 \leq 1, \quad
x_{\frac{n}{2}+1}^4+\cdots+x_{n}^4 \leq 1,
\end{align*}
where $n$ is even. Since the degree is $6$, we apply the  $3^{rd}$
Lasserre's relaxation \reff{consos:deg2d}. The cone $\mc{K}$ has $3$
blocks. For $n=20$, solving
\reff{prmsdp:consos}-\reff{dualsdp:consos} takes about $3.2$ hours.
The computed solution of the SDP relaxation has error around $1
\cdot 10^{-7}$. The computed lower bound $f_{sos}^{con} \approx
-38.8840$. The computed optimal $Z^*$ has rank one, and we get a
global optimal solution.
Its error is around $1 \cdot 10^{-7}$.
\qed \end{exm}

\begin{exm}
Minimize the quartic polynomial
\[
\sum_{1\leq i < j \leq n} \left( x_ix_j+x_i^2x_j-x_j^3-x_i^2x_j^2
\right)
\]
over the hypercube $[-1,1]^n=\{x \in \re^n:  x_i^2 \leq 1\}$. We
apply the $2^{nd}$ Lasserre's relaxation \reff{consos:deg2d}. The
resulting cone $\mc{K}$ has $n+1$ blocks. For $n=50$, solving
\reff{prmsdp:consos}-\reff{dualsdp:consos} takes about $2.8$ hours.
The error of the computed solution for the SDP relaxation is around
$10^{-6}$. The computed lower bound $f_{sos}^{con} \approx -1250$.
The computed optimal $Z^*$ does not satisfy FEC.  So, we use
$Z^{*}(2:n+1,1)$ as a starting point,
and get a local minimizer (by function {\tt fmincon} in Matlab Optimization
Toolbox)
\[
\baray{l}
(-1,-1,-1,-1,-1,-1,-1,-1,-1,-1,-1,-1,-1,1,1,1,-1,-1,-1,1,1, \\
1,-1,1,-1,1,-1,-1,1,1,1,1,-1,1,-1,1,1,1,1,1,1,1,1,1,1,1,1,1,1,1).
\earay
\]
Its objective value is $f(x^{*}) = -1232$, which is greater than the
lower bound $f_{sos}^{con}$. This is probably because the second
order Lasserre's relaxation itself is not exact.
\qed \end{exm}

\begin{exm}[Random polynomials]
We test the performance of Algorithm~\ref{alg:Newton-CG}
(implemented by {\tt SDPNAL}) in minimizing polynomials over the
unit ball. Generate $f(x)$ randomly as
\[
f(x) = \sum_{\af\in \N^n: |\af|\leq 2d }  f_\af x^\af,
\]
where the coefficients $f_\af$ are Gaussian random variables. We
solve the SDP relaxation \reff{prmsdp:consos}-\reff{dualsdp:consos}
by {\tt SDPNAL}. The cases of degrees $4,6,8,10$ are tested. The
computational results are in Table~\ref{constab:rd:dg=4:10}.
\begin{table}[htb]
\centering
\begin{tiny}
\begin{tabular}{|c||c||r|r|r||c||c|} \hline
($n$,$2d$)  &  \#Inst& \multicolumn{3}{c||}{time (min, med, max)} &
errsol(min, med, max)&  errsdp(min, med, max) \\ \hline (30,4)& 10 &
0:00:28   & 0:00:52 & 0:02:47   & (5.6e-8, 1.3e-6, 6.9e-6) &
(1.3e-7, 8.1e-7, 2.9e-6)
\\ \hline (40,4) &
10 &  0:03:35 & 0:06:38  & 0:10:32  & (8.8e-8, 1.8e-6, 9.5e-6)
&(2.2e-7, 1.0e-6, 4.5e-6)
\\ \hline (50,4) &    3 &  0:20:34 & 0:22:18  &  0:24:59 &
  (5.7e-6, 5.6e-6, 7.0e-6)  & (2.7e-6, 2.8e-6, 3.4e-6)  \\ \hline
  (60,4) &    3 &  0:35:02 & 1:20:15  &  1:20:38 &
  (1.5e-7, 3.5e-6, 2.5e-5)  & (1.7e-7, 1.7e-6, 1.2e-5)  \\ \hline
(20,6) &   3 & 0:36:31  & 0:49:17  & 0:56:35 &(8.5e-7, 2.7e-6,
4.4e-6)    & (5.8e-7, 1.3e-6, 2.7e-6)\\ \hline
 (12,8) &  3 &  0:27:11  & 0:44:06 & 0:59:30   & (5.5e-7, 2.8e-6, 9.0e-6)  & (9.0e-7, 1.3e-6, 4.2e-6) \\
\hline (9,10) &  3 &  0:16:31  & 0:36:05 & 0:40:53   & (2.6e-7,
3.3e-6, 1.4e-5)  & (2.7e-7, 1.6e-6, 6.3e-6) \\ \hline
  (80,4) &   3 &  10:52:30 & 15:12:40  &  15:57:30 &
 (5.3e-6, 5.5e-6, 2.2e-1)  & (2.6e-6, 2.6e-6, 2.7e-3)  \\ \hline
 (25,6) &  3 & 10:38:04  & 11:00:48  &
 12:57:59 &(5.9e-3, 6.6e-3, 1.4e-2)   & (3.6e-3, 5.8e-3, 6.1e-3)  \\
\hline
\end{tabular}
\end{tiny} \caption{Computational results for random constrained polynomial
optimization} \label{constab:rd:dg=4:10}
\end{table}
When $(n,2d)=(80,4)$ or $(25,6)$, the SDP relaxations are not solved
very well sometimes. This is probably because of the incurred
ill-conditioning. In all the other cases, the SDP relaxations are
solved quite well, and accurate global minimizers are found.
%
\qed \end{exm}

\section{Some discussions}
\setcounter{equation}{0}

In this section, we discuss some numerical issues about the
performance of regularization methods in solving SDP relaxations for
large scale polynomial optimization problems.

\subsection{Scaling polynomial optimization}

SDP relaxations arising from polynomial optimization are harder to
solve than general SDP problems. A reason for this is that the
polynomials are not scaled very well sometimes. For instance, if the
optimal $Z^{*} $ has rank 1, then $Z^{*} = [x^{*}]_d[x^{*}]_d^T$
($x^*$ is a minimizer) has entries of the form
\[
1, x^{*}_{1},\ldots,(x^{*}_{1})^{2},\ldots,
\ldots,(x^{*}_{1})^{2d},\ldots,(x^{*}_{n})^{2d}.
\]
Clearly, if some coordinate $x_i^*$ is small or big, then $Z^*$ is
badly scaled and its entries $Z^*_{ij}$ easily get
underflow/overflow during the computation. This might cause severe
ill-conditioning in computations and make the computed solutions
less accurate.
Scaling is a useful approach to overcome this issue. In
\cite{Gloptipoly,PS01}, it was pointed out that scaling is important
in solving polynomial optimization problems efficiently. Generally, there is
no simple rule to select the best scaling factor. In the following, we
propose a practical scaling procedure.

%
Let $s=(s_1,\ldots, s_n)>0$ and scale $x$ to $\hat{x} = (\hat{x}_1,
\ldots, \hat{x}_n)$ as
\[
x \quad = \quad (s_1\hat{x}_1, \ldots, s_n\hat{x}_n).
\]
Then $f(x)$ is scaled to be the polynomial $f(s_1\hat{x}_1, \ldots,
s_n\hat{x}_n)$ in $\hat{x}$. The best scaling factor $s$ should be
such that the global minimizers of the scaled polynomial have
coordinates close to one or negative one. This is difficult because
optimizers are usually unknown.
However, as the algorithm runs, one often gets close to minimizers
and would estimate them from the computations. Typically, we only
need to scale the problem when the algorithm fails to converge.
Sometimes, we might need to do scaling several times. From our
experiences, a practical scaling procedure is:

\bdes

\item [Step 1]  If Algorithm \ref{alg:Newton-CG} converges well, we do no
scaling and let it run; otherwise, select a scaling
vector $s=(s_1,\ldots, s_n)>0$ as:

\be \label{factor} s_{i}  = \bca
\tau  & \text{ if }   |y_{e_i}| \leq \tau, \\
|y_{e_i}| & \text{ otherwise.} \eca \ee Here $\tau>0$ is fixed and
$y$ is the most recent update for an optimal $y^*$ of
\reff{sdp:momsos}.

%

\item [Step 2] Scale $f(x)$ as
$f(s_1\hat{x}_1,\ldots, s_n\hat{x}_n)$. Go back to Step 1 and solve
the scaled polynomial optimization again.


\edes

In the above, $\tau>0$ is usually (but not too) small, because the
coefficients of the scaled polynomial $f(s_1\hat{x}_1, \ldots,
s_n\hat{x}_n)$ should not be very tiny. We use $\tau = 10^{-3}$ in
the examples below.

\begin{exm} \label{exm:5.1}
Consider the polynomial optimization \be \label{P:S:em}
\min_{x\in\re^n} \quad  x^{4}_{1} + \ldots + x^{4}_{n}+
\sum\limits_{1\leq i<j<k\leq n}x_{i}x_{j}x_{k}. \ee
For this kind of polynomials, its global minimizers usually have
large negative values and lead to ill-conditioning of the SDP
relaxation (cf. \cite[\S 5.1]{PS01}). Here we show the importance of
scaling for the case $n=20$.
\begin{table}[htb]
\centering

  \begin{tabular}{|c|c|c|c|} \hline
  Iter & time &  low. bdd. &  sdp err.  \\
  \hline 1 &    0:01:15   &  -1.0806e+7 &    0.7555 \\
  \hline 2 &
  0:01:15   &  -1.9444e+7 &   0.0460 \\
 \hline 3 &    0:00:37   &
 -2.1883e+7 &   0.0082 \\
 \hline 4 &    0:01:16   &  -2.2266e+7 &
 2.4e-6 \\ \hline
 \end{tabular}
  \caption{Results of scaling process for Example~\ref{exm:5.1}.}
 \label{stumfel:parillo}
 \end{table}


We use the scaling procedure described above. The computational
results are shown in Table \ref{stumfel:parillo}. The ``low. bdd."
there stands for the computed optimal value of SDP relaxation
\reff{sdp:stdsos}-\reff{sdp:momsos}, which is always a lower bound
of the global minimum, and ``sdp err." stands for the error of the
computed solution of \reff{sdp:stdsos}-\reff{sdp:momsos}, which is
defined in \reff{def:sdp-err}.
It takes four times of scaling to solve the SDP relaxation
reasonably well.
\qed
\end{exm}

\begin{exm} \label{exm:5.2}
Consider the least square problem (Watson function \cite{MGH}): \be
\label{nonlinear:least:squre} \min_{x\in \re^{n}} \quad
\sum\limits_{i=1}^{m} f^{2}_{i}(x) . \ee Here $n=30$ and the
polynomials $f_i$ are defined as follows:
\begin{eqnarray}
 \begin{array}{l}
f_{i}(x) = \sum\limits_{j=2}^{n}(j-1)x_j t_i^{j-2}-
\left(\sum\limits_{j=1}^{n}x_j t^{j-1}_i \right)^2-1,\,\, t_i =
\frac{i}{29},~~ 1\leq i \leq 29,
\end{array}
\end{eqnarray}
and $f_{30} = x_1,~~ f_{31} = x_2 -x^{2}_{1} -1$.
Its SDP relaxation \reff{sdp:stdsos}-\reff{sdp:momsos} has size
$(N,m)=(496,46376)$. We solve it by the scaling procedure mentioned
earlier.
\begin{table}[htb]
\centering

 \begin{tabular}{|c|c|c|c|} \hline
  Iter & time &  low. bdd. &   sdp err.  \\ \hline 1 &  0:28:09 &
  -9.1556 &    0.9955 \\ \hline 2&    0:45:34  &    0.0134  & 8.1e-3
  \\ \hline 3&    0:30:40   &     0.1468  &   9.1e-4 \\ \hline 4&
  0:26:33    &    0.1298 &    4.7e-4 \\ \hline 5&   0:25:18    &
 0.0969  &   3.1e-4 \\ \hline 6&   0:18:10   &   0.0648  &    8.3e-5
 \\ \hline
 \end{tabular}
 \caption{Results of scaling process for Example~\ref{exm:5.2}. }
  \label{watson:scaling:iteration}
 \end{table}


The results are in Table~\ref{watson:scaling:iteration}. It takes
six times of scaling to solve the SDP relaxation reasonably well.
\qed \end{exm}

\begin{remark}
In each step of the scaling process,
we need to solve a new SDP problem of the same size as the earlier
one. As shown in Section 4, sometimes the SDP relaxations in
polynomial optimization could be solved very well without scaling.
But this is not always the case, e.g., like Examples
\ref{exm:5.1} and \ref{exm:5.2}.
Typically, we need to do scaling only when Algorithm~\ref{alg:Newton-CG} has troubles to solve
a problem. The performance of Algorithm~\ref{alg:Newton-CG} is bad when the
SDP problem is ill-conditioned or has degeneracy. Our experiments
show that scaling can help solve the problem more efficiently.
\end{remark}

\subsection{Why do we need regularization methods
in large scale polynomial optimization?}

%

As we have seen in Introduction, a major issue of interior point
methods is that in each step one needs to solve an $m\times m$
linear system and two $N\times N$ matrix equations. This would be a
big restriction in applications if $m$ is huge, because it requires
storing an $m\times m$ matrix in computer and $\mc{O}(m^3)$
arithmetic operations. Unfortunately, SDP relaxations from
polynomial optimization have an unfavorable property that
$m=\mc{O}(N^2)$. As shown in Table~\ref{tab:(N,m)}, in minimizing a
general quartic polynomial of $100$ variables, the SDP relaxation
has $m$ greater than $4$ million. To solve such an SDP relaxation by
interior point methods, one needs to store a square matrix of length
bigger than $4$ million in memory. On a regular computer, this is
almost impossible. However, regularization methods requires much
less memory storage. In Algorithm~\ref{alg:Newton-CG}, in each inner
loop, we still need to solve the linear system \reff{linsys:Newton}
which is also $m\times m$. But, the Hessian $\nabla_y^2
\varphi_{\sigma}(Y,y)$ does not need to be explicitly formulated.
Actually, the authors of \cite{ZST08} showed that the matrix vector
product $\nabla_y^2 \varphi_{\sigma}(Y,y) \cdot z$ would be
evaluated in $\mc{O}(m)$ arithmetic operations and the memory
requirement for \reff{linsys:Newton} has linear order in $m$.
Because of this special feature, CG type methods are very suitable
for solving \reff{linsys:Newton}. This property has been
successfully used by software {\tt SDPNAL}.

For unconstrained polynomial optimization, its SDP relaxation has an
attractive feature. From the construction of $A_\af$ in
subsection~\ref{sec:unopt}, we can easily see that distinct
$A_\af$'s have no common nonzero entries. Thus, the matrices $A_\af$
in \reff{def:ucpop-A} are orthogonal to each other, and the matrix
$\mc{A}\mc{A}^*$ is diagonal. In this case, Algorithm~\ref{BPM:Ksdp}
is easily implementable, because its every step only involves
solving a diagonal linear system and computing an eigenvalue
decomposition. In Algorithm \ref{alg:Newton-CG}, the diagonal
$(\mc{A}\mc{A}^*)^{-1}$ would be used as a preconditioner for
\reff{linsys:Newton} in CG iterations.

\subsection{Other numerical methods}

To the authors' best knowledge, there are few efficient numerical methods for
solving large scale polynomial optimization problems. One method
that might be useful in applications is the low rank method proposed by Burer and
Monteiro \cite{Bour03} (implemented in software {\tt SDPLR}
\cite{SDPLR}). In some cases, the dual optimal $Z^*$ of SDP
relaxations might have low rank. Thus, in such situations, {\tt
SDPLR} would be applied to solve the dual SDP relaxation like
\reff{sdp:momsos} or \reff{dualsdp:consos} (not the primal SDP
relaxation \reff{sdp:stdsos} or \reff{prmsdp:consos}, since $X^*$
typically has high rank). We tested {\tt SDPLR} on some
examples in this paper. Its performance is similar to {\tt SDPNAL}.
However, {\tt SDPLR} is less attractive theoretically and suitable
only when $Z^*$ has low rank. This is because the basic idea of {\tt
SDPLR} is to change SDP into a nonlinear programming problem via
matrix factorization, and typically one would only get a local
optimizer. However, by {\tt SDPLR}, it is not guaranteed
to get an optimizer of the SDP relaxation. Moreover,
even if an optimizer of SDP is obtained, its
optimality can not be certified. A reason for this is that {\tt
SDPLR} is not a primal-dual type method, and typically a primal-dual
pair is required to check optimality.
On the other hand, the computational performance of {\tt SDPLR} is promising.
It is an interesting future work to investigate
properties of the low rank method in solving polynomial optimization.

There are interesting recent work on solving
large scale polynomial optimization problems by other methods.
Bertsimas, Freund and Sun \cite{DFS11} proposed an accelerated first order method
to solve unconstrained polynomial optimization problems.
A nice theoretical property of first order type methods is that
there are bounds on the complexity of computations, as proved in \cite{DFS11}.
Henrion and Malick \cite{HM112,HM11} proposed a projection method
for solving conic optimization and SOS relaxations.
These methods can solve bigger problems than the interior point
methods do, but they might take a big number of iterations to get
an accurate optimal solution and generally its convergence is slow.
In practical computations of solving big polynomial optimization problems,
Algorithm~\ref{alg:Newton-CG} typically has faster convergence,
because it uses second order information
(e.g., approximate Newton directions).

\subsection{Convergence and nondegeneracy}

The performance of Algorithm~\ref{alg:Newton-CG}
is not always very good for solving SDP relaxations in
polynomial optimization. As we have seen earlier, a typical reason
is the ill-conditioning. Another reason might be the degeneracy of
the SDP relaxations. In \cite{ZST08}, it was shown that if the SDP
problem is nondegenerate, then Algorithm~\ref{alg:Newton-CG} has
good convergence; otherwise, it might converge very badly or even
does not converge. Generally, it is difficult to check in advance
whether an SDP relaxation is degenerate or not. For SDP relaxation
\reff{sdp:stdsos}-\reff{sdp:momsos} in unconstrained polynomial optimization,
or \reff{prmsdp:consos}-\reff{dualsdp:consos} in constrained
optimization, a typical case for it to be degenerate is that a
polynomial optimization problem has several distinct global minimizers. To
see this for the unconstrained polynomial optimization \reff{uc:minf(x)},
suppose it has two distinct global minimizers $u^*,v^*$ and the SOS
relaxation \reff{opt:stdsos} is exact. Then, the optimal values of
\reff{sdp:stdsos} and \reff{sdp:momsos} are equal, and
\reff{sdp:momsos} has two distinct optimal $Z^*$ (being
$[u^*]_d[u^*]_d^T$ and $[v^*]_d[v^*]_d^T$). This implies the primal
SDP relaxation \reff{sdp:stdsos} is degenerate. The situation is
similar for constrained polynomial optimization. From this
observation, Algorithm~\ref{alg:Newton-CG} might not be very efficient if
the SDP relaxation is exact and there are more than one distinct
optimizers. Of course, Algorithm~\ref{alg:Newton-CG} might still
work if an SDP problem is degenerate, like in Example~4.1. But this
is occasional and typically not the case in practice.

\bigskip
\bigskip
\noindent {\bf Acknowledgement} \quad The authors would like to
thank Xinyuan Zhao, Defeng Sun and
Kim-Chuan Toh for sharing their software {\tt SDPNAL}, and Gabor
Pataki for comments on the degeneracy of SDP. They also thank Bill
Helton and Igor Klep for fruitful discussions on this work.

\end{document}